\documentclass [12] {article}
\usepackage{amsthm}
\usepackage{epsfig, graphicx, xcolor}
\usepackage{amssymb}
\usepackage{amsmath, amsfonts}
\newtheorem{theorem}{Theorem}

\addtocounter{page}{0}

\begin{document}

\title{
Approximation of some classes of set-valued periodic functions by generalized trigonometric polynomials}

\author{V. F.~Babenko$^*$, V. V.~Babenko$^{**}$, M. V.~Polischuk$^*$}

\date{$^*$Oles Honchar Dnepropetrovsk National University,\\
$^{**}$The University of Utah$\;\;\;\;\;\;\;\;\;\;\;\;\;\;\;\;\;\;\;\;\;\;\;\;\;\;\;\;\;\;\;\;\;\;\;\;\;\;\;\;\;\;\;\;\;\;$\\
{\it E-mail: babenko.vladislav@gmail.com},\\
{\it E-mail: vera.babenko@gmail.com},\\
        {\it E-mail: polishchuk.mariya@gmail.com}}

\maketitle

{\large
\sloppy
\begin{abstract}



Generalizations of some known results on the best, best linear and best one-sided approximations by trigonometric polynomials of the classes of $2\pi$ - periodic functions presented in the form of convolutions to the case of set-valued functions are obtained

\bigskip


\end{abstract}

\section{Introduction}

In approximation theory, results on the exact solutions of problems of the best, best linear and best one-sided approximation, on the classes of periodic functions by trigonometric polynomials, are well-known. Review and details of most of the results in this area, as well as further references, can be found in the following articles ~\cite{Nik},~\cite{Dz} and monographs~\cite{Korneichuk},~\cite{KLD}. The purpose of this work is to generalize some of the results to the case of set-valued functions.

Problems of approximation of set-valued functions have been studied relatively recently. Overview and some known results can be found in ~\cite{Vit},~~\cite{Artstein},~\cite{DynF},~\cite{DynFM}.


Let us briefly describe the structure of the work. In the second section of the article we will provide some necessary definitions, notations and facts that are related to the case of numerical-valued periodic functions. In the third section we will present necessary definitions and facts from the theory of set-valued functions. In the section four we will give formulations of the problems of approximation of the set-valued functions. Fifth section is dedicated to some results of approximation of set-valued periodic functions that can be represented as a convolution, by some $"$generalized trigonometric polynomials$"$.

\section{Approximation of the classes of numerical functions}

Let $C$ and $L_p\,(1\le p\le \infty)$ be the spaces of $2\pi$ - periodic functions ${f: \mathbb{R}\to \mathbb{R}}$ endowed with the corresponding norms $\|\cdot\|_C$  and $\|\cdot\|_{Lp}$.

Let $X$ be $L_p,\, 1\le p <\infty ,$ or $C$, and let $H$ be finite-dimensional subspace of the space $X$. For $f\in X$ let
\begin{equation}
\label{funkK}
E\left(f,H\right)_{X}=\inf_{{T}\in H}\left\| f-{T}\right\|_{X}.
\end{equation}
Let also ${\cal M}\subset X$ be some class of functions,
\begin{equation}
\label{funkF}
{E}\left({\cal M},H\right)_{X}=\sup_{f\in {\cal M}} E_n\left(f,H\right)_X.
\end{equation}
Values (\ref{funkK}) and (\ref{funkF}) are called the best approximation of the function $f$ and of the class ${\cal M}$ respectively, by subspace $H$ in the space $X$.

Furthermore, let
$$
{U}\left({\cal M},H\right)_{X}=\inf\limits_{A}\sup_{f\in {\cal M}} \| f-Af\|_X,
$$
where $\inf_{A}$ is taken over all possible linear operators $A:X\to H$. The value ${U}\left({\cal M},H\right)_{X}$ is called the best linear approximation of the class ${\cal M}$ by subspace $H$ in the space $X$.

As usual, the convolution $K*\varphi$ of the functions $K\in L_1 $ (kernel of the convolution) and $\varphi\in L_1$, we will define by the following equality
\begin{equation}
\nonumber
K*\varphi(x)= \int\limits_0^{2\pi}K(t)\varphi(x-t)dt.
\end{equation}
Let $F_p=\{\varphi\in L_p,1\le p\le \infty: \|\varphi\|_{p}\le 1\}$. Denote by $K*F_p$ the class of functions that have the following form
\begin{equation}
\nonumber
f(x)=K*\varphi(x),\quad \varphi\in F_p.
\end{equation}
It is well known (see, e.g. ~\cite{Nik},~\cite{Dz},~\cite{Korneichuk}~\cite{Bab_smzh},~\cite{BL}), that many important classes of real-valued periodic functions are classes of the form $K*F_p$.

Denote by $H^{{T}}_{2n-1}$, $n=1,2\ldots$, the set of all trigonometric polynomials ${T}_{n-1}(x)$ of order $n-1$, i.e. the set of functions of the form
$$
T_{n-1}(x)=\frac{a_0}{2}+\sum\limits_{k=0}^{n-1}a_k\cos kx+b_k\sin kx,\; a_k, b_k\in\mathbb{R}.
$$

S. M. Nikolskii~\cite{Nik} gave rather general conditions on a kernel  $K$ that allows to find the values ${E}\left(K*F_\infty,H^{{T}}_{2n-1}\right)_{C}$ and ${E}\left(K*F_1,H^{{T}}_{2n-1}\right)_{L_1}$.
We present the condition $N_n^*$ only:

{\it We will say that a kernel $K$ satisfies the condition  $N_n^*$ if there exist a polynomial $T^*\in H^T_{2n-1}$ and a point $\theta\in [0,\pi /n]$ such that almost everywhere
$$(K(x)-T^*(x))\varphi_n(x-\theta)\ge 0.$$
}
Here and everywhere below
$$
\varphi_{n}(x):={\rm sgn}\sin nx.
$$

S. M. Nikolskii have proved the following

\noindent{\bf Theorem A.} {\it If a kernel $K$ satisfies the condition $N_n^*$, then
$$
{E}\left(K*F_\infty,H^{{T}}_{2n-1}\right)_{C}={U}\left(K*F_\infty,H^{{T}}_{2n-1}\right)_{C}
$$
$$
={E}_n\left(K*F_1,H^{{T}}_{2n-1}\right)_{L_1}={U}_n\left(K*F_1,H^{{T}}_{2n-1}\right)_{L_1}
$$
$$
={E}\left(K,H^{{T}}_{2n-1}\right)_{L_1}=\| K-T^*\|_{L_1}=\| K*\varphi_{n}\|_{C}.
$$
}

Many kernels, that are important for the approximation theory, satisfy the condition $N_n^*$ (one can find corresponding examples in  ~\cite{Nik},~\cite{Dz},~\cite{Bab_smzh},~\cite{BL}).

\bigskip

 Regarding best linear approximations of the classes of periodic functions, we present here the following theorem that can be proved as the theorem 1 from~\cite{B_Pich}.

\bigskip

\noindent{\bf Theorem B.} {\it Let $p,q\in [1,\infty]$ and $p^{-1}+q^{-1}=1$. If $K\in L_q$, then for any $n\in\mathbb{N}$
$$
U(K*F_p,H^T_{2n-1})_C=E(K,H^T_{2n-1})_{L_q}.
$$
}

\section{Necessary definitions and facts related to the set-valued functions}

Let us present some necessary definitions and facts related to the space of nonempty, compact subsets of the space $\mathbb{R}^{m}$, and  set-valued functions.
Proofs of these facts can be found in в~\cite{Bal},~\cite{aubin-fra}~\cite{HuPap}.
Denote by $K(\mathbb{R}^{m} )$ the space of nonempty, compact subsets of the space $\mathbb{R}^{m}$. By $K^c(\mathbb{R}^{m} )$ we will denote the set of all convex elements of the space $K(\mathbb{R}^{m} )$.
We will consider set-valued $2\pi$ --periodic functions $f:\mathbb{R}\to K(\mathbb{R}^{m} )$ i.e., functions  such that $f(x+2\pi)=f(x)$ for any  $x\in\mathbb{R}$.

As usual, a linear combination of sets $A,B\subset K(\mathbb{R}^{m} )$ is defined by
\[\lambda A+\mu B=\left\{\lambda a+\mu b:a\in A,b\in B\right\},
\quad\lambda ,\mu \in \mathbb{R}.
\]
Convex hull of the set $A\subset
K(\mathbb{R}^{m} )$ we will denote by $\mathrm{co}A$.

If $a=(a_1,...,a_m)\in \mathbb{R}^m$ and $\xi=(\xi_1,\ldots,\xi_m)\in \mathbb{R}^m$,  then
$$|a|_{l^m_2}:=\sqrt{\sum\limits_{j=1}^{m}{ a_{j}}^2 },\;\;(a,\xi ):=\sum\limits_{k=1}^m a_k\xi_k.$$ For a point $a\in
\mathbb{R}^m$, and some set $B\in K(\mathbb{R}^{m} )$, define
$$d(a,B):=\mathop{\inf }\limits_{b\in B} |a-b|_{l^m_2}.$$ This is the distance from the point $a$ to the set  $B$. For sets $A,B\in
K(\mathbb{R}^{m} )$ let
$$d(A,B):=\mathop{\sup }\limits_{a\in A} d(a,B).$$ This is the distance from the set $A$ to the set $B$. Hausdorff metric $\delta$ in the space $K(\mathbb{R}^{m} )$ is defined in the following way. If $A, B\in K(\mathbb{R}^{m} )$, then \[\delta(A,B):=\max \{ d(A,B),d(B,A)\}. \]
Note, that $K(\mathbb{R}^{m} )$ and $K^c(\mathbb{R}^{m} )$ with Hausdorff metric are complete metric spaces.

Metric $\delta\left(A,B\right)$ has following properties
\[
 \delta\left(\lambda A,\lambda B\right)=\lambda \delta\left(A,B\right),  \quad\quad  \forall \lambda >0,\quad \forall A, B\in K(\mathbb{R}^{m}
 ),\]
\[\delta\left(A+B,C+D\right)\le \delta\left(A,C\right)+\delta\left(B,D\right), \quad \quad \forall A,\, B,\, C,\, D\in K(\mathbb{R}^{m}
),\]
\[\delta (\mathrm{co}\; A,\mathrm{co}\; B) \le \delta(A,B), \quad\quad \forall A, B\in K(\mathbb{R}^{m} ),
\]
\[
\delta (\alpha A, \beta A)\le |\alpha -\beta |\|A\|,\quad\quad \forall \alpha , \beta\in \mathbb{R},\forall A\in K^c(\mathbb{R}^m),
\]
where $\| A\|:=\delta (A,\{ \theta\})$ (here and below $\theta =(0,\ldots ,0)$ is the null element of the space $\mathbb{R}^m$).

\bigskip

The Aumann's integral of a set-valued function $f:[0,2\pi ]\to K({\mathbb R}^m)$ is defined as the set of all integrals of  integrable selections of $f$:
\[
I(f)=\int\limits_0^{2\pi}f(x)dx:=\left\{ \int\limits_0^{2\pi}\phi (x)dx\; :\;
\phi (x)\in f(x)\; {\rm a.\, e.},\; \phi \; {\rm is \;
integrable}\right\}.
\]
It's known (see, e.g,~\cite{Aseev}) that if the function $f:[0,2\pi ]\to K({\mathbb R}^m)$ is measurable and the function  $\| f(\cdot)\|$ is integrable (set of such functions we will denote by $L^A_1$), then
\[
\int\limits_0^{2\pi}f(x)dx\in K^c({\mathbb R}^m).
\]
We will need following properties of Aumann's integral for functions $f, g\in L^A_1$ (see, e.g.,~\cite{aubin-fra},~\cite{HuPap}):
\[\displaystyle\int\limits_0^{2\pi}\mathrm{co} (f(x))dx=\int\limits_0^{2\pi}f(x)dx,\]
\[\displaystyle\delta\left(\int\limits_ 0^{2\pi}f(x)dx,\int\limits_ 0^{2\pi}g(x)dx\right)\leq \int\limits_ 0^{2\pi}\delta(f(x),g(x))dx.
\]

\section{Problems of approximation that are related to set-valued functions}
Denote by $L^A_p,\; 1\le p\le \infty,$ the set of functions $f\in L^A_1$, such that $\| f(\cdot )\|\in L_p$. In $L^A_p$ introduce a metric by setting
\[
\delta_{L^A_p}(f,g):= \| \delta (f(\cdot),g(\cdot))\|_{L_p}.
\]
Let also
$$
\Phi_p :=\{ f\in L^A_p\; :\; \delta_{L^A_p}(f(\cdot),\{\theta \})\le 1\}.
$$

Let a kernel $K\in L_1$ be given. We will consider problems of approximation of classes  $K*\Phi_p$ of set-valued functions, that can be written in the form
$$
f(x)=\int\limits_0^{2\pi}K(x-t)g(t)dt,\;\; g\in \Phi_p,
$$
where the integral is considered in the Aumann's sense. Due to the properties of the Aumann's integral, functions from the class $K*\Phi_p$ are convex-valued.

As approximating functions, we will use set-valued functions of the form
\begin{equation}
\label{sv_trig}
\tau (x)=\int\limits_0^{2\pi}T(x-t)h(t)dt,
\end{equation}
where $T$ is a polynomial from $H^T_{2n-1}$, and $h$ is a function from $L^A_1$. The set of all possible functions of such form we will denote by $SVH^T_{2n-1}$.

For the function $f\in K*\Phi_p, \; 1\le p\le \infty,$ define
\begin{equation}
\nonumber
\mathcal{E}\left(f,SVH^T_{2n-1}\right)_{L^A_p}=\inf\limits_{\tau\in SVH^T_{2n-1}} \delta_{L^A_p}\left(f,\tau\right).
\end{equation}
Let also
\begin{equation}
\nonumber
\mathcal{E}\left(K*\Phi_q,SVH^T_{2n-1}\right)_{L^A_p}=\sup_{f\in K*\Phi_q}\mathcal{E}\left(f,SVH^T_{2n-1}\right)_{L^A_p}.
\end{equation}

The problem of finding the value $\mathcal{E}\left(K*\Phi_p,SVH^T_{2n-1}\right)_{L^A_p}$ is a set-valued analog of the problem of best approximation of the class of functions by trigonometric polynomials.

We will also consider a problem of finding following values
$$
\mathcal{U}\left(K*\Phi_p,SVH^T_{2n-1}\right)_{L^A_q}:=\inf\limits_{T\in H^T_{2n-1}}\sup\limits_{f=K*g\in K*\Phi_p}\| K*g-T*g\|_{L_q^A}.
$$
This problem is a set-valued analog of the problem of best linear approximation of corresponding classes of numerical-valued functions.

Now we extend the set of approximating functions and instead of the set $SVH^T_{2n-1}$ we will use set $\widetilde{SVH}^T_{2n-1}$ of set-valued functions of the form
\begin{equation}\label{sv_trig+}
\tilde{\tau} (x)=\tau (x)+B_r(\theta),
\end{equation}
where $\tau \in SVH^T_{2n-1}$ and $B_r(\theta)=\{ z\in\mathbb{R}\; :\; |z|_{l^m_2}\le1\},\; r\ge 0$.

For function $f\in K*\Phi_p, \; 1\le p\le \infty,$ define
\begin{equation}
\nonumber
\mathcal{E}^+\left(f,\widetilde{SVH}^T_{2n-1}\right)_{L^A_p}:=\inf\limits_{\tau\in \widetilde{SVH}^T_{2n-1}\atop \forall x\; f(x)\subset \tau (x)} \delta_{L^A_p}\left(f,\tau\right).
\end{equation}
The problem of finding (estimating) of the value
\begin{equation}
\nonumber
\mathcal{E}^+\left(K*\Phi_q,\widetilde{SVH}^T_{2n-1}\right)_{L^A_p}=\sup_{f\in K*\Phi_q}\mathcal{E}^+\left(f,\widetilde{SVH}^T_{2n-1}\right)_{L^A_p}
\end{equation}
is a set-valued analog of the problem of best one-sided approximation on the class of numerical-valued functions and on our opinion has some interest.

\section{Results}
\begin{theorem}
Let $n\in \mathbb{N}$, $p,q\in [1,\infty]$ and $p^{-1}+q^{-1}=1$.
If $K\in L_q$, then
\begin{equation}
\label{otsenka1}
\mathcal{E}\left(K*\Phi_p,SVH^{T}_{2n-1}\right)_{L_\infty^A}\le \mathcal{U}\left(K*\Phi_p,SVH^T_{2n-1}\right)_{L^A_\infty}\le  E_n\left(K,H^{{T}}_{2n-1}\right)_{L_q}.
\end{equation}
If $K\in L_1$, then for $p\in [1,\infty]$
\begin{equation}
\label{otsenka}
\mathcal{E}\left(K*\Phi_p,SVH^{T}_{2n-1}\right)_{L_p^A}\le \mathcal{U}\left(K*\Phi_p,SVH^T_{2n-1}\right)_{L^A_p}\le  E_n\left(K,H^{{T}}_{2n-1}\right)_{L_1}.
\end{equation}

\end{theorem}

$\Box$
Let us first prove  the inequality (\ref{otsenka1}).
Let $T^*\in H^T_{2n-1}$ be the polynomial of the best $L_q$-approximation for $K$ and let $g\in \Phi_p$. Using the properties of the Hausdorff metric and Aumann's integral we estimate $\delta (K*g(x),T^**g(x))$:
$$
\delta (K*g(x),T^**g(x))=\delta\left(\int\limits_0^{2\pi}K(x-t)\; g(t)dt,\int\limits_0^{2\pi}T^*(x-t)\;g(t)dt\right)
$$
$$
=\delta\left(\int\limits_0^{2\pi}K(x-t)\; {\rm co}\,g(t)dt,\int\limits_0^{2\pi}T^*(x-t)\;{\rm co}\,g(t)dt\right)
$$
$$
\le \int\limits_0^{2\pi}\delta\left(K(x-t){\rm co}\,g(t),T^*(x-t){\rm co}\,g(t)\right)dt
$$
$$
\le\int\limits_0^{2\pi}|K(x-t)-T^*(x-t)|\delta\left({\rm co}\,g(t),\{\theta\}\right)dt.
$$
Therefore,
\begin{equation}
\label{upper}
\delta (K*g(x),T^**g(x))\le \int\limits_0^{2\pi}|K(x-t)-T^*(x-t)|\delta\left({\rm co}\,g(t),\{\theta\}\right)dt.
\end{equation}
Using H\"{o}lder inequality we obtain
$$
\|\delta (K*g(\cdot),T^**g(\cdot))\|_{L_\infty}
$$
$$
\le \max\limits_{x\in\mathbb{R}}\left(\int\limits_0^{2\pi}|K(x-t)-T^*(x-t)|^qdt\right)^{\frac 1q}\left(\int\limits_0^{2\pi}\delta\left({\rm co}\,g(t),\{\theta\}\right)^pdt\right)^{\frac 1p}\le
E_n\left(K,H^{{T}}_{2n-1}\right)_{L_q}.
$$
We have proved the relation (\ref{otsenka1}).

Let now $K\in L_1$, $f=K*g\in K*\Phi_p$, and let $T^*\in H^T_{2n-1}$ be the polynomial of the best $L_1$-approximation for $K$.
Using the inequality (\ref{upper}) and generalized Minkovskii inequality we obtain
$$
\| \delta (K*g(\cdot),T^**g(\cdot))\|_{L_p}\le
\left\| \int\limits_0^{2\pi}|K(\cdot -t)-T^*(\cdot -t)|\delta\left({\rm co}\,g(t),\{\theta\}\right)dt\right\|_{L_p}
$$
$$
\le \| K-T^*\|_{L_1}\| g(\cdot )\|_{L^A_p}\le E(K,H^T_{2n-1})_{L_1}.
$$
Therefore,
$$
\mathcal{E}\left(K*\Phi_p,SVH^{T}_{2n-1}\right)_{L_p^A}\le \mathcal{U}\left(K*\Phi_p,SVH^T_{2n-1}\right)_{L^A_p}\le  E_n\left(K,H^{{T}}_{2n-1}\right)_{L_1}.
$$
Relation (\ref{otsenka}) is also proved. Thus Theorem 1 is proved.
$\Box$

\begin{theorem}
If a kernel $K(t)$ satisfies the condition $N_n^*$, then for $p=1$ or $p=\infty$ the following relations hold
\begin{equation}\nonumber
\mathcal{E}\left(K*\Phi_p,SVH^{T}_{2n-1}\right)_{L^A_p}=\mathcal{U}\left(K*\Phi_p,SVH^T_{2n-1}\right)_{L^A_p}
$$
$$
={E}\left(K*F_p,H^{{T}}_{2n-1}\right)_{L_p}=\|K*\varphi_n\|_{L_\infty},
\end{equation}
where
$\varphi_n(x)={\rm sign} \sin n x$, $x\in {\mathbb R}$.

\end{theorem}

$\Box$
For $p=1$ and $p=\infty$ the estimate
$$
\mathcal{E}\left(K*\Phi_p,SVH^{T}_{2n-1}\right)_{L^A_p}\le \mathcal{U}\left(K*\Phi_p,SVH^T_{2n-1}\right)_{L^A_p}
$$
$$
\le E\left(K,H^{{T}}_{2n-1}\right)_{L_1}=\|K*\varphi_n\|_{L_\infty}
$$
follows from theorem 1 and theorem A.

We will first obtain  the estimate from below in the case $p=\infty$.
Choose arbitrary $a\in \mathbb{R}^{m}$, such that $\delta^h(\{ a\},\{\theta\})=|a|=1$.
Then $K*(\varphi_n(\cdot)\{ a\})=K*\varphi_n(\cdot)\cdot\{ a\}\in\Phi_\infty$.

Let also, $\tau\in SVH^T_{2n-1}$ has the form (\ref{sv_trig}) and let $\psi$ be arbitrary selection from the function $h$. Then
$$
\| \delta(K*\varphi_n(\cdot)\cdot\{ a\},\tau(\cdot))\|_{L_\infty}=\max_{x\in\mathbb{R}} \delta(K*\varphi_n(x)\cdot\{ a\},\tau(x))
$$
$$
\ge \max_{x\in\mathbb{R}} d(\tau(x),K*\varphi_n(x)\cdot\{ a\})
\max_{x\in\mathbb{R}}d\left( \int\limits_{0}^{2\pi}T(x-t)h(t)dt,K*\varphi_n(x)\cdot\{ a\}\right)
$$
$$
\ge\max_{x\in\mathbb{R}}\left|\int\limits_{0}^{2\pi}T(x-t)\psi(t)dt-K*\varphi_n(x)\cdot\{ a\}\right|_{l^m_2}
$$
$$
= \max_{x\in\mathbb{R}}\sup_{\xi\in\mathbb{R}^m\atop |\xi |_{l_2^m}=1}\left|\left(\xi,\int\limits_{0}^{2\pi}T(x-t)\psi(t)dt\right)-K*\varphi_n(x)\cdot (\xi, a)\right|
$$
$$
= \max_{x\in\mathbb{R}}\sup_{\xi\in\mathbb{R}^m\atop |\xi |_{l_2^m}=1}\left|\int\limits_{0}^{2\pi}T(x-t)\left(\xi,\psi(t)\right)dt-K*\varphi_n(x)\cdot (\xi, a)\right|
$$
$$
\ge \sup_{\xi\in\mathbb{R}^m\atop |\xi |_{l_2^m}=1}|(\xi ,a)|\max_{x\in\mathbb{R}}\left| K*\varphi_n(x)\right|=\| K*\varphi_n\|_{L_\infty}
$$
(the last inequality holds due to the Chebyshev's alternance theorem). Consequently,
$$
\mathcal{E}\left(K*\Phi_\infty ,SVH^{T}_{2n-1}\right)_{L^A_\infty}=\sup_{f\in K*\Phi_\infty}\inf\limits_{{\tau}\in SVH^{T}_{2n-1}} \|\delta\left(f(\cdot),{\tau}(\cdot)\right)\|_\infty
$$
$$
\ge \inf\limits_{{\tau}\in SVH^{T}_{2n-1}}\| \delta(K*\varphi_n(\cdot)\cdot\{ a\},\tau(\cdot))\|_{L_\infty}\ge \| K*\varphi_n\|_{L_\infty}.
$$
Estimate from below, as well as, the statement of the theorem for the case $p=\infty$ are proved.

Obtain now an estimate from below for $p=1$. Choose arbitrary $a\in \mathbb{R}^{m}$, such that  $\delta(\{ a\},\{\theta\})=|a|=1$, and arbitrary function $g\in F_1$. It is clear that $g(\cdot)\cdot \{ a\}\in \Phi_1$ and $f(\cdot)=K*g(\cdot)\cdot\{ a\}\in K*\Phi_1$.

Let $\tau\in SVH^T_{2n-1}$ be of the form (\ref{sv_trig}) and let $\psi$ be arbitrary integrable selection from the function $h$. We have
$$
\|\delta (f(\cdot),\tau(\cdot))\|_{L_1}=\int\limits_0^{2\pi}\delta (K*g(x)\cdot \{ a\},\tau(x))dx
$$
$$
\ge \int\limits_0^{2\pi}d (\tau(x), K*g(x)\cdot \{ a\})dx =
\int\limits_0^{2\pi} d\left( \int\limits_{0}^{2\pi}T(x-t)h(t)dt,K*g(x)\cdot\{ a\}\right)dx
$$
$$
\ge\int\limits_0^{2\pi} \left| \int\limits_{0}^{2\pi}T(x-t)\psi(t)dt-K*g(x)\cdot\{ a\}\right|_{l^m_2}dx
$$
$$
=\int\limits_0^{2\pi}\sup_{\xi\in\mathbb{R}^m\atop |\xi |_{l^m_2}=1}  \left|\left(\xi, \int\limits_{0}^{2\pi}T(x-t)\psi(t)dt\right)-K*g(x)\cdot(\xi, a)\right|dx
$$
$$
=\int\limits_0^{2\pi}\sup_{\xi\in\mathbb{R}^m\atop |\xi |_{l^m_2}=1}  \left| \int\limits_{0}^{2\pi}T(x-t)\left(\xi,\psi(t)\right)dt-K*g(x)\cdot(\xi, a)\right|dx
$$
$$
\ge\sup_{\xi\in\mathbb{R}^m\atop |\xi |_{l^m_2}=1}  \int\limits_0^{2\pi}\left| \int\limits_{0}^{2\pi}T(x-t)\left(\xi,\psi(t)\right)dt-K*g(x)\cdot(\xi, a)\right|dx
$$
$$
\ge \sup_{\xi\in\mathbb{R}^m\atop |\xi |=1}|(\xi ,a)|E(K*g, H^T_{2n-1})_{L_1}=E(K*g, H^T_{2n-1})_{L_1}.
$$

Therefore,
$$
\mathcal{E}\left(K*g\cdot \{ a\} ,SVH^{T}_{2n-1}\right)_{L^A_1}\ge E\left(K*g ,H^{T}_{2n-1}\right)_{L_1},
$$
and taking into account theorem A we obtain
$$
\mathcal{E}\left(K*\Phi_1, SVH^{T}_{2n-1}\right)_{L^A_1}\ge E\left(K*F_1 ,H^{T}_{2n-1}\right)_{L_1}=\| K*\varphi_n\|_{L_\infty}.
$$
We obtain estimate from below in the case $p=1$. Theorem is proved.
$\Box$

Estimate from below relating to the case $p=1$ allows the following generalization.

\begin{theorem}
For any kernel $K\in L_1$, any $p,q\in [1,\infty],\, q\le p,$ and any $n\in \mathbb{N}$
\begin{equation}\nonumber
\mathcal{E}\left(K*\Phi_p,SVH^{T}_{2n-1}\right)_{L^A_q}\ge {E}\left(K*F_p,H^{{T}}_{2n-1}\right)_{L_q}.
\end{equation}
\end{theorem}

$\Box$
Choose arbitrary $a\in \mathbb{R}^{m}$ such that $|a|_{l_2^m}=1$ and arbitrary function $g\in F_p$. It is clear that $g(\cdot)\cdot \{ a\}\in \Phi_p$ and $f(\cdot)=K*g(\cdot)\cdot\{ a\}\in K*\Phi_p$.

As in the proof of the previous theorem, let $\tau\in SVH^T_{2n-1}$ be of the form (\ref{sv_trig}) and let $\psi$ be an arbitrary integrable selection from $h$. We have
$$
\|\delta (f(\cdot),\tau(\cdot))\|_{L_q}=\left(\int\limits_0^{2\pi}\delta (K*g(x)\cdot \{ a\},\tau(x))^qdx\right)^{\frac 1q}
$$
$$
\ge \left(\int\limits_0^{2\pi}d (\tau(x), K*g(x)\cdot \{ a\})^qdx \right)^{\frac 1q}
$$
$$
=\left(\int\limits_0^{2\pi}d\left( \int\limits_{0}^{2\pi}T(x-t)h(t)dt,K*g(x)\cdot\{ a\}\right)^qdx\right)^{\frac 1q}
$$
$$
\ge\left(\int\limits_0^{2\pi} \left| \int\limits_{0}^{2\pi}T(x-t)\psi(t)dt-K*g(x)\cdot\{ a\}\right|_{l^m_2}^qdx\right)^{\frac 1q}
$$
$$
=\left(\int\limits_0^{2\pi}\sup_{\xi\in\mathbb{R}^m\atop |\xi |_{l^m_2}=1}  \left|\left(\xi, \int\limits_{0}^{2\pi}T(x-t)\psi(t)dt\right)-K*g(x)\cdot(\xi, a)\right|^q dx\right)^{\frac 1q}
$$
$$
\ge\sup_{\xi\in\mathbb{R}^m\atop |\xi |_{l^m_2}=1} \left(\int\limits_0^{2\pi} \left|\int\limits_{0}^{2\pi}T(x-t)\left(\xi, \psi(t)\right) dt-K*g(x)\cdot(\xi, a)\right|^q dx\right)^{\frac 1q}
$$
$$
\ge \sup_{\xi\in\mathbb{R}^m\atop |\xi |=1}|(\xi ,a)|E(K*g, H^T_{2n-1})_{L_q}=E(K*g, H^T_{2n-1})_{L_q}.
$$

Thus
$$
\mathcal{E}\left(K*g\cdot \{ a\} ,SVH^{T}_{2n-1}\right)_{L^A_q}\ge E\left(K*g ,H^{T}_{2n-1}\right)_{L_q},
$$
and, therefore,
$$
\mathcal{E}\left(K*\Phi_p ,SVH^{T}_{2n-1}\right)_{L^A_q}\ge E\left(K*F_p ,H^{T}_{2n-1}\right)_{L_q}.
$$
Theorem is proved.
$\Box$

The following theorem gives us set-valued analog of the theorem B.
\begin{theorem}
Let $n\in \mathbb{N}$, $p,q\in [1,\infty]$ and $p^{-1}+q^{-1}=1$.
If $K\in L_q$, then
\begin{equation}
\label{otsenka3}
\mathcal{U}\left(K*\Phi_p,SVH^T_{2n-1}\right)_{L^A_\infty}=  E_n\left(K,H^{{T}}_{2n-1}\right)_{L_q}.
\end{equation}
\end{theorem}

$\Box$
The inequality
$$
\mathcal{U}\left(K*\Phi_p,SVH^T_{2n-1}\right)_{L^A_\infty}\le  E_n\left(K,H^{{T}}_{2n-1}\right)_{L_q}
$$
follows from the theorem 1. Let us prove the inequality of the opposite sense.

Choose arbitrary $a\in \mathbb{R}^{m}$, such that  $|a|_{l_2^m}=1$. It is clear that for any function  $g\in F_p$, we will have $g(\cdot)\cdot \{ a\}\in \Phi_p$ and  $f(\cdot)=K*g(\cdot)\cdot\{ a\}\in K*\Phi_p$.

For any $T\in H^T_{2n-1}$ we will have
$$
\delta (K*g(x)\{ a\},T*g(x)\{ a\})=| K*g(x)-T*g(x)|\cdot |a|_{l^m_2}=| K*g(x)-T*g(x)|.
$$
Therefore,
$$
\sup\limits_{g\in \Phi_p}\|\delta (K*g(x),T*g(x))\|_{L_\infty}\ge \sup\limits_{g\in F_p}\|\delta (K*g(x)\{ a\},T*g(x)\{ a\})\|_{L_\infty}
$$
$$
=\sup\limits_{g\in F_p}\| K*g(x)-T*g(x)\|_{L_\infty}=\sup\limits_{g\in F_p}\max\limits_{x\in\mathbb{R}}\left|\int\limits_0^{2\pi}(K(x-t)-T(x-t)g(t)dt\right|
$$
$$
=\max\limits_{x\in\mathbb{R}}\sup\limits_{g\in F_p}\left|\int\limits_0^{2\pi}(K(x-t)-T(x-t))g(t)dt\right|=\max\limits_{x\in\mathbb{R}}\| K(x-\cdot)-T(x-\cdot)\|_{L_q}
$$
$$
=\| K-T\|_{L_q}\ge E(K,H^T_{2n-1})_{L_q}.
$$
Theorem is proved.
$\Box$

In conclusion we present a theorem which gives an estimate from above for the value $\mathcal{E}^+\left(K*\Phi_p,\widetilde{SVH}^{T}_{2n-1}\right)_{L_\infty^A}$.
\begin{theorem}
Let $n\in \mathbb{N}$, $p,q\in [1,\infty]$ and $p^{-1}+q^{-1}=1$.
If $K\in L_q$, then
\begin{equation}
\nonumber
\mathcal{E}^+\left(K*\Phi_p,\widetilde{SVH}^{T}_{2n-1}\right)_{L_\infty^A}\le 2 E_n\left(K,H^{{T}}_{2n-1}\right)_{L_q}.
\end{equation}
\end{theorem}

$\Box$
Let $T^*\in H^T_{2n-1}$ be a polynomial of the best $L_q$-approximation for $K$. Let also, a function $f=K*\Phi_p$ be given. Due to the theorem 1 we have for any $x\in\mathbb{R}$
\begin{equation}\nonumber
\delta (K*g(x),T^**g(x))\le E(K,H^T_{2n-1})_{L_q}=:e.
\end{equation}
Set
$$
\widetilde{\tau}(x)=T^**g(x)+B_e(\theta).
$$
It is easily seen that $$K*g(x)\subset \widetilde{\tau}(x)$$ and
$$
\delta (K*g(x),\widetilde{\tau}(x))=\delta (K*g(x),T^**g(x)+B_e(\theta))
$$
$$
\le\delta (K*g(x),T^**g(x))+\delta (B_e(\theta),\theta)\le 2E(K,H^T_{2n-1})_{L_q}.
$$
Theorem 5 follows from this estimate.
$\Box$

}


\begin{thebibliography}{99}

\bibitem {Nik} {Nikolskii S.~M.} Approximation of functions in the mean by trigonometrical polynomials// Izv. AN SSSR, Ser. mat., 1946, v. 10,  N 3, p. 207--256.
\bibitem{Dz} {Dzyadyk V.~K.} On best approximation in classes of periodic functions defined by integrals of a linear combination of absolutely monotonic kernels//Mat. Zametki, 1974, v.~16, N~5, p. 691--701.

\bibitem{Korneichuk} Korneichuk N.P., Exact constants in approximation theory. M.: Nauka, 1987, 424 p.
\bibitem{KLD} Korneichuk N. P., Ligun A. A., Doronin V. G., Approximation with constraints.-- Kiev, Naukova dumka, 1982, 252 p.
\bibitem{Bab_smzh} Babenko V. F.,  Approximation of convolution classes //Sib. Mat. Zhurn, 1987, v. 28, N 5, p. 6 - 21.

\bibitem{BL} Babenko V. F., Ligun A.A.,  Development of studies on the exact solution of extremal problems of the theory of best approximation //Ukr. Mat. Zhurn, 1990, v. 42, N 1, p. 4 - 17.

\bibitem{B_Pich} Babenko V. F., Pichugov S. A., Best linear approximation of some classes of differentiable periodic functions//Mat. Zametki, 1980, v. 27, N 5, p. 683 - 689.

\bibitem{Vit}

  {Vitale R.~A.} Approximations of convex set-valued functions// J. Approxim. Theory, 1979, 26, P.~301--316.

\bibitem{Artstein}
  { Artstein Z.}
Piecewise linear approximations of set-valued maps//J. Approxim. Theory, 1989,  56, P.~41--47.

\bibitem{DynF}  { Dyn N., Farkhi  E.,}
Approximations of set-valued functions with compact images -- an overview, approximation and probability// Banach Center Publ., 2006, Vol.~72, P.~1--14.

\bibitem{DynFM}
Dyn, N., Farkhi, E., Mokhov, A.  Approximation of
set-valued functions: Adaptation of classical approximation
operators. Hackensack: Imperial College Press, 2014, 153 p.


\bibitem{Bal} {Polovinkin Е. S., Balashov М. V.} Elements of convex and strongly convex analysis. M.: Fizmatlit, 2004. -- 416 p.


\bibitem{aubin-fra} Aubin, J. P., Frankowska, H. Set-valued analysis.
Modern Birkhauser Classics: Boston, 1990, 461 p.

\bibitem{HuPap}

Hu S., Papageorgiou N., Handbook of multivalued analysis, Vol. 1: Theory, Kluwer Academic Publishers, 1997, 964 p.


\bibitem{Aumann} {Aumann R.J.} Integrals of set-valued functions// J. Math. Anal. and Appl., 1965, Vol. 12, N 1, P.1-12.



\bibitem{Aseev} {Aseev S. M.} Quasilinear operators and their application in the theory of multivalued mappings// Trudy MIAN SSSR.- 1985, v. 167, p. 25 – 52.



\end{thebibliography}
\end{document}